\documentstyle[leqno,amssymb,times]{article}

\def\a{\alpha}
\def\b{\beta}

\long\def\th#1#2#3{\vskip-\lastskip\vskip4pt plus2pt
\noindent{\uppercase{#1}} #2\hskip-\lastskip\ {\it #3}\vskip-\lastskip\vskip4pt plus2pt}
\def\Proof{\vskip-\lastskip\vskip4pt plus2pt
{\noindent\uppercase{Proof}.\ }\ignorespaces}
\def\endproof{\unskip\nobreak\kern5pt\nobreak\vrule height4pt width4pt depth0pt
\vskip4pt plus2pt}
\def\sq{\unskip\nobreak\kern5pt\nobreak\vrule height4pt width4pt depth0pt}

\def\CU{{\mathcal U}}
\def\CA{{\mathcal A}}
\def\CH{{\mathcal H}}
\def\cp{\Delta}
\def\oh{\frac{1}{2}}
\def\ooh{\frac{3}{2}}
\def\sqq{\sqrt{q}}
\def\isq{\frac{1}{\sqq}}
\def\ts{\otimes}
\def\oq{\frac{1}{q}}
\def\acts{\triangleright}
\def\ket#1{|#1\rangle}
\def\endproof{\vrule height 0.5em depth 0.2em width 0.5em}
\newbox\tbox
\newbox\aubox
\newbox\adbox
\def\title#1{\setbox\tbox=\hbox{\let\\=\cr
\baselineskip14pt\vbox{\Large\bf\tabskip 0pt plus15cc
\halign to\hsize{\hfil\ignorespaces \uppercase{##}\hfil\cr#1\cr}}}}
\newbox\abbox
\setbox\abbox=\vbox{\vglue18pt}
\def\author#1{\setbox\aubox=\hbox{\let\\=\cr
\baselineskip12pt\vbox{\tabskip 0pt plus15cc
\halign to\hsize{\hfil\ignorespaces \uppercase{{##}}\hfil\cr#1\cr}}}%
\global\setbox\abbox=\vbox{\unvbox\abbox\box\aubox\vskip8pt}}
\def\address#1{\setbox\adbox=\hbox{\let\\=\cr
\baselineskip12pt\vbox{\it\tabskip 0pt plus15cc
\halign to\hsize{\hfil\ignorespaces {##}\hfil\cr#1\cr}}}%
\global\setbox\abbox=\vbox{\unvbox\abbox\box\adbox\vskip16pt}}
\def\makemytitle{
\begin{center}
\box\tbox
\box\abbox
\end{center}}

\begin{document}
\title{Dirac Operator on the Standard \\ Podle\'s Quantum Sphere}
\author{Ludwik D\c{a}browski}
\address{Scuola Internazionale Superiore di Studi Avanzati,\\
Via Beirut 2-4, I-34014, Trieste, Italy}
\author{Andrzej Sitarz}
\address{Institute of Physics, Jagiellonian University,\\
Reymonta 4, 30-059 Krak\'ow, Poland}
\makemytitle
\abstract{
Using principles of quantum symmetries we derive the
algebraic part of the real spectral triple data for
the standard Podle\'s quantum sphere: equivariant representation,
chiral grading $\gamma$, reality structure $J$ and the
Dirac operator $D$, which has bounded commutators with
the elements of the algebra and satisfies the first order
condition.}

\vspace{4mm}
{\small
\section*{1. Introduction.}

It is a prevailing common opinion that Connes' approach
to noncommutative geometry [4] is in a sense incompatible
with the (compact) quantum-group examples of $q$-deformations [9].
The main obstacle was the construction of an appropriate
Dirac operator $D$, fulfilling the requirements postulated in [5].
In this paper we solve this problem for the standard Podle\'s
quantum sphere [14] by finding a Dirac operator satisfying
the `algebraic' part of the axioms of real spectral triples.
We simplify our considerations by making an additional assumption
of the symmetry with respect to the $\CU_q(su(2))$ action, thus
pursuing the principles of construction of equivariant real
spectral triples [15].

\section*{2. The standard quantum sphere and its symmetry.}

We recall the definition of the algebra $\CA(S^2_q)$
of the standard Podle\'s quantum sphere as a subalgebra
of the quantum $SU(2)$ function algebra [14,16].

Let $q$ be a real number $0<q \leq 1$ and $\CA(SU_q(2))$ be a
$*$-algebra generated by $a$, $a^*$ and $b,b^*$, which
satisfy the following commutation rules:

\begin{equation}
\begin{array}{llll}
& ba = q ab,  &&  \\
& b^*a = qab^*, && bb^*=b^*b, \\
& a^*b = q ba^*, && a^*a + q^2 bb^* = 1, \\
& a^*b^*=q b^* a^*, && aa^* + bb^*=1.\\
\end{array}
\end{equation}
The algebra $\CA(S^2_q)$ is isomorphic to the subalgebra
generated by:
\begin{equation}
\begin{array}{l}
B=  a b, \\
B^* = b^* a^*, \\
A = bb^*.
\end{array}
\end{equation}
They obey the following relations:
\begin{equation}
\begin{array}{lll}
\label{srel}
AB = q^2 BA, & \phantom{xxx} & A B^* = q^{-2} B^* A, \\
B B^* = q^{-2}A (1 - A), & & B^* B = A (1 - q^2 A).
\end{array} \label{PodSta}
\end{equation}
The quantized algebra $\CU_q(su(2))$ has $e,f,k,k^{-1}$ as
generators of the $*$-Hopf algebra, satisfying relations:
\begin{equation}
\begin{array}{lllll}
ek = qke,  &\phantom{xxx}& kf = qfk, & \phantom{xxx}
& k^2 - k^{-2} = (q-q^{-1})(fe-ef).\end{array}
\end{equation}
The coproduct is given by:
\begin{equation}
\begin{array}{lllll}
\cp k = k \ts k,  &\phantom{xxx}&\cp e =  e \ts k + k^{-1} \ts e,
&\phantom{xxx}& \cp f =  f \ts k + k^{-1} \ts f .
\end{array}
\end{equation}
The counit $\epsilon$, antipode $S$, and star structure are as folows:
\begin{equation}
\begin{array}{lllll}
 \epsilon(k) = 1,  &\phantom{xxx}&  \epsilon(e) = 0,
&\phantom{xxx}& \epsilon(f) = 0, \\
Sk = k^{-1}, &\phantom{xxx}& Sf = - qf, &\phantom{xxx}& Se = -q^{-1} e, \\
 k^* = k, &\phantom{xxx}& e^* =f, &\phantom{xxx}& f^* =e.
\end{array}
\end{equation}
{}From the usual Hopf algebra pairing between $\CU_q(su(2))$ and $\CA(SU_q(2))$
we obtain an action of $\CU_q(su(2))$ on $\CA(SU_q(2))$ given on generators by:
\begin{equation}
\begin{array}{lll}
\label{Uact}
 k \acts a = q^{\frac{1}{2}} a,  &\phantom{xxx}&    k \acts a^* = q^{-\frac{1}{2}} a^*,\\
 k \acts b = q^{\frac{1}{2}} b, &\phantom{xxx}& k \acts b^* = q^{-\frac{1}{2}} b^*,\\
 e \acts a = -b^*, &\phantom{xxx}& e \acts a^* = 0, \\
 e \acts b = q^{-1} a^*, &\phantom{xxx}& e \acts b^* = 0, \\
 f \acts a = 0, &\phantom{xxx}& f \acts a^* = q b, \\
 f \acts b = 0, &\phantom{xxx}& f \acts b^* = -a.
\end{array}
\end{equation}
This action preserves the $*$-structure:
\begin{equation}
h \acts (x^*) = \left( (Sh)^* \acts x \right)^*, \;\
\forall h \in  \CU_q(su(2)),
x \in \CA(S^2_q). \label{staract}
\end{equation}
Hence we derive the action of $\CU_q(su(2))$ on the generators of the
standard Podle\'s sphere:
\begin{equation}
\begin{array}{l}
e \acts B= - (q^{\frac{1}{2}} + q^{-\frac{3}{2}}) A + q^{-\frac{3}{2}},
\\
e \acts B^* = 0,
\\
e \acts A = q^{-\frac{1}{2}} B^*, \\
k \acts B = q B,
\\
k \acts B^* = q^{-1} B^*,
\\
k \acts A = A, \\
f \acts B =0,
\\
f \acts B^* = (q^{\frac{3}{2}} + q^{-\frac{1}{2}}) A - q^{-\frac{1}{2}},  \\
f \acts A = - q^{\frac{1}{2}} B.
\end{array} \label{aphereact}
\end{equation}

\section*{3. Equivariant representation of \boldmath{$\CA(S^2_q)$}.}

As the next step we shall find a representation of the algebra
$\CA(S^2_q)$ on a Hilbert space $\CH$, which is equivariant
under the action of $\CU_q(su(2))$ defined
in (\ref{aphereact}).
\th{Definition 1.}{}{
Let $V$ be an $A$-module and $H$ be a Hopf algebra.
Also, let $V$ and $A$ be $H$-modules.
We call $V$ an $H$-equivariant $A$-module
if the following condition is satisfied
\begin{equation}
\label{covar}
 h (\alpha v) = (h_{(1)} \acts \alpha ) (h_{(2)} v)\ ,
~~~\forall h \in H, ~\alpha \in A, ~v\in V.
\end{equation}
Here we use Sweedler's notation for the coproduct of $H$
and $\acts$ for the action of $H$ on $A$.
(In other words, $V$ is a module over the crossed (smash) product of
$H$ and $A$.)
}
\th{Definition 2.}{}{
A bounded representation $\pi$ of $A$ on a Hilbert space $\CH$
is called $H$-equivariant if there exists a dense linear subspace
$V$ of $\CH$ such that $V$ is an  $H$-equivariant $A$-module
and $\pi (\a ) v = \a v, ~~\forall
v\in V, ~\alpha \in A$.
}
\vspace{3mm}
We shall use the (known) representation theory of $\CU_q(su(2))$.
\th{Lemma 3.}{({[9]})}{
The irreducible finite dimensional representations of
$\CU_q(su(2))$ are labeled by \mbox{$l=0,\frac{1}{2},1,\ldots$}
and they are given by
\begin{equation}
\begin{array}{l}
 f \ket{l,m} = \sqrt{[l-m][l+m+1]} \ket{l,m+1}\\
 e \ket{l,m} = \sqrt{[l-m+1][l+m]} \ket{l,m-1}\\
 k \ket{l,m}  = q^{m} \ket{l,m}
\end{array}
\end{equation}
where $m\in\{-l, -l\! +\! 1, \dots ,l\! -\! 1, l\}$ and for any number $x$
$$ [x] :=\frac{q^x-q^{-x}}{q-q^{-1}}.$$
For a given $l$ we denote such space $V_l$.
}
\noindent
The representations above are $*$-repre\-senta\-tions of
$\CU_q(su(2))$ with respect to the (hermitian) scalar product for which
the vectors $\ket{l,m}$ are orthonormal.

\vspace{3mm}
We attempt now to construct a $\CU_q(su(2))$-equivariant
representation of $\CA(S^2_q)$ on the Hilbert space completion
of a $\CU_q(su(2))$-equivariant $\CA(S^2_q)$-module constructed
as a direct sum of all representations $V_l$.
A priori some multiplicities can occur but we make a simplifying assumption
that this is not the case, guided by the picture of spinors
on the classical (commutative) sphere.
However, it may be easily verified that the results
of Lemma 3 could be generalized to the case with
multiplicities.
\vspace{2mm}
\th{Lemma 4.}{}{
A\,  $\CU_q(su(2))$-equivariant representation
of $\CA(S^2_q)$ on the Hilbert space completion of
$V_0\nobreak \oplus V_\frac{1}{2} \oplus V_1 \oplus \ldots$~
must have the following form:
\begin{eqnarray}
\label{recurs}
\label{rep1}
 B \ket{l,m} = B^+_{l,m} \ket{l+1,m+1} +B^0_{l,m}
\ket{l,m+1} + B^-_{l,m} \ket{l-1,m+1}, \label{eqdefB} &&\\
\label{rep2}
 B^*\ket{l,m} =  \widetilde{B}^+_{l,m} \ket{l+1,m-1}
+ \widetilde{B}^0_{l,m} \ket{l,m-1} +
\widetilde{B}^-_{l,m}\ket{l-1,m-1}, \label{eqdefC}&&\\
\label{rep3}
 A \ket{l,m}  =  A^+_{l,m} \ket{l+1,m} + A^0_{l,m}
\ket{l,m}  +  A^-_{l,m} \ket{l-1,m}, \label{eqdefA}&&
\end{eqnarray}
where $A^j_{l,m},B^j_{l,m},\tilde{B}^j_{l,m}$ for $j=+,0,-$,
are constants.}
\Proof We use the covariance property (\ref{covar})
on the algebraic direct sum of $V_j$.
First, from the action of $k$ on $B$ we see that $B$ must increase
the exponent $m$ of the eigenvalue of $k$ by $1$.
Similar arguments show that $B^*$
decreases $m$ by $1$ and $A$ does not change $m$.
Therefore $B \ket{l,m}$ is a sum:
\begin{equation}
B \ket{l,m} = \sum_i B_i \ket{i,m+1}, \label{dum1}
\end{equation}
where the sum runs over $i = 0, \oh, 1, \ooh, \dots$
with finite number of nonzero coefficients.

Applying $f^{l-m+1}$ to both sides we
obtain $0$ on the left-hand side and a sum starting with
$i=l+2$ on the right-hand side.  By linear independence
of all elements in this sum, we see that the sum over
$i$ in (\ref{dum1}) runs only up to $l+1$.
(This does not eliminate yet the possibility of $i<l-1$.)

Using a similar argument, with $B^*$ instead of $B$, we observe
that again the terms with $i>l+1$, in a corresponding sum, do not appear.
Since we are looking for a $*$-representation of $\CA(S^2_q)$,
we conclude that in both expressions $i$ runs from $l-1$ to $l+1$, as
stated in (\ref{eqdefB}-\ref{eqdefA}). From the sphere defining
relations (\ref{PodSta}) it follows that this holds also for $A$.
\sq

It follows from (\ref{rep1}-\ref{rep3}) that we can separate the half integer spin
representations from the integer spin representations.
In the sequel we restrict ourselves only to the case of the Hilbert space
$\CH_{\oh}$
given by completion of the direct sum
$$V = V_{\oh} \oplus V_{\ooh} \oplus \dots $$
of all half-integer representations of $\CU_q(su(2))$,
motivated by the classical picture of (chiral) spinors over $S^2$.

We proceed by applying $f$ on both sides of (\ref{eqdefB}) to obtain
the following recursion relations for the coefficients $A^j_{l,m},B^j_{l,m}$:
\begin{equation}
\label{recurss}
\frac{1}{q} B^j_{l,m+1} \sqrt{[l-m][l+m+1]} =
B^j_{l,m} \sqrt{[l+j-m-1][l+j+m+2]}.
\end{equation}
We solve them explicitly for each value of $m$:
\begin{equation}
\begin{array}{l}
\phantom{xxxx} B^+_{l,m}  = q^{m} \sqrt{ [l+m+1][l+m+2]} \;\alpha^+_l, \\
\phantom{xxxx} B^0_{l,m}  = q^{m} \sqrt{ [l+m+1][l-m]} \;\alpha^0_l, \\
\phantom{xxxx} B^-_{l,m}  = q^{m} \sqrt{ [l-m][l-m-1]} \;\alpha^-_l,
\end{array}
\label{Brep}
\end{equation}
where $\alpha^\pm_l, \alpha^0_l$ are yet undetermined functions of $l$.

If we apply $e$ to (\ref{eqdefC}), we obtain
\begin{equation}
\begin{array}{l}
\phantom{xxxx} \widetilde{B}^+_{l,m}  = q^{m-1} \sqrt{ [l-m+2][l-m+1]} \; \alpha^-_{l+1},  \\
\phantom{xxxx} \widetilde{B}^0_{l,m}  = q^{m-1} \sqrt{ [l+m][l-m+1]} \; \alpha^0_l, \\
\phantom{xxxx} \widetilde{B}^-_{l,m}  = q^{m-1} \sqrt{ [l+m][l+m-1]}\; \alpha^+_{l-1}.
\end{array}
\label{Crep}
\end{equation}

Finally, for $A$, using the action of $e$ on $B$ (or equivalently
$f$ on $B^*$) we get:
\begin{equation}
\begin{array}{l}
\phantom{xxxx} A^+_{l,m}  = -q^{m+l+\frac{1}{2}} \sqrt{ [l-m+1][l+m+1]} \;\alpha^+_l \\
\phantom{xxxx} A^0_{l,m}  = q^{-\oh} \frac{1}{1+q^2} \left( [l-m+1][l+m]  - q^2 [l-m][l+m+1]
\right) \alpha^0_l + \frac{1}{1+q^2}, \\
\phantom{xxxx}  A^-_{l,m}  =  q^{m-l-\frac{1}{2}} \sqrt{ [l-m][l+m]} \;\alpha^-_l \ .
\end{array} \label{Arep}
\end{equation}
Moreover, the $*$-representation condition yields a consistency relation
between $\alpha^+$ and
$\alpha^-$:
\begin{equation}
\alpha^-_{l+1} = - q^{2l+2} \alpha^+_l\ ,  ~~{\rm for}~~ l = \oh, \ooh, \dots \ .
\label{ab1}
\end{equation}

Now, we use the relations defining the algebra $\CA(S^2_q)$
in order to calculate the form of $\alpha^+_l$ and $\alpha^0_l$.
The comparison of the $\ket{l+1,m+1}$ component of the action
of the first relation in (\ref{srel}) on
$\ket{l,m}$ gives the following recurrence relation in $l$:
\begin{equation}
\alpha^0_{l+1} [2l+4] = \alpha^0_l [2l] + \isq (q - q^{-1}).
\label{0recur}
\end{equation}

To get the relation for $\alpha^+_l$ and obtain the initial
terms (that is values for $l=\frac{1}{2}$) we
employ the last two relations in (\ref{srel}), of which we choose the following
convenient combination:
\begin{equation}
\label{combi}
B^* B - q^4 B B^* = (1 - q^2) A. \label{rel1}
\end{equation}

Our results can be summarized in the following lemma.
\th{Lemma 5.}{}{
There are two inequivalent $\CU_q(su(2))$-equivariant representations
of $\CA(S^2_q)$ on the Hilbert space $\CH_\oh $
given by formulae (\ref{eqdefB}-\ref{eqdefA}),
(\ref{Brep}-\ref{ab1}) and by the following formulae
for $\alpha^+_l$ and $\alpha^0_l$:
\begin{itemize}
\item $\pi_+$:
\end{itemize}
\begin{eqnarray}
\alpha^0_l = \isq \frac{ (q - \oq)
[l - \oh][l+ \frac{3}{2}] + q}{[2l][2l+2]}\ , \\
\alpha^+_l = q^{-l-2}
\frac{1}{\sqrt{[2l+2]([4l+4]+[2][2l+2])}}\ .
\end{eqnarray}
~\\
\begin{itemize}
\item $\pi_-$:
\end{itemize}
\begin{eqnarray}
\alpha^0_l = \isq \frac{ (q - \oq)
[l - \oh][l+ \frac{3}{2}] - q^{-1}}{[2l][2l+2]}\ , \label{a01}\\
\alpha^+_l = q^{-l-1} \frac{1}{\sqrt{[2l+2]([4l+4]+[2][2l+2])}}\ .
\label{ap1}
\end{eqnarray}
}
\Proof By direct calculation with the help of the symbolic
algebra program.\footnote{The source file with definitions and
examples is available from the web-site: \\
\ \  \hbox to 6.5mm{}
{\footnotesize\tt http://www.cyf-kr.edu.pl/\~{}ufsitarz/maple.html}}

First we verify the value of
$\alpha^+_{\frac{1}{2}}$:
\begin{equation}
( \alpha^+_{\frac{1}{2}} )^2 = q^{-4 \mp 1} \frac{1}{([3])^2 [4]},
\end{equation}
where the $\pm$ signs corresponds to the freedom we have. This
leads to the explicit solution of the recurrence relation
(\ref{0recur}), which is the formula (\ref{a01}). For $\alpha^+_l$
we use two relations obtained from (\ref{rel1}):
\begin{equation}
\begin{array}{l}
(\alpha^+_l)^2 q^{2l+3} [2l+3] [2]
- (\alpha^+_{l-1})^2 q^{2l+1} [2l-1] [2] \\
~\\
- (\alpha^0_l)^2 q [4l+2]/[2l+1]
+ \alpha^0_l \sqrt{q}\, (q - q^{-1}) = 0,
\end{array}
\end{equation}
~\\
\begin{equation}
\begin{array}{l}
-(\alpha^+_l)^2 q^{2l+3} [4l+6] [2]
+(\alpha^+_{l-1})^2 q^{2l+1} [4l-2] [2] \\
~\\
+ (\alpha^0_l)^2 q ([2])^2
- \alpha^0_l \sqrt{q}\, (q - q^{-1}) [4l+2]/[2l+1]
+ (q- q^{-1})^2 = 0.
\end{array}
\end{equation}
By subtracting the left-hand sides to eliminate $\alpha^+_{l-1}$
and using the result (\ref{a01}) for $\alpha^0_l$,
we arrive at the final formula (\ref{ap1}) presented above.
Besides (\ref{combi}) we also check that the last two relations
in (\ref{srel}) are separately fulfilled
and that all elements of $\CA(S^2_q)$ are
represented as bounded operators on $\CH_\oh $.

Note that the possible sign ambiguity in the formula for $\alpha^+_l$
can be globally resolved by the redefinition of the basis. \sq

\section*{4. Real spectral triple for \boldmath{$S^2_q$}.}

We begin by taking as the Hilbert space of Dirac spinors,
on which $\CA(S^2_q)$ is represented,
the direct sum $\CH := \CH_{\oh} \oplus \CH_{\oh}$,
with the representation given by:
\begin{equation}
\pi(a) = \left( \begin{array}{ll} \pi_+(a) & 0 \\ 0 & \pi_-(a)
\end{array} \right),
\end{equation}
and the grading $\gamma$:
$$
\gamma = \left( \begin{array}{ll} 1 & 0 \\ 0 & -1 \end{array} \right).
$$
It is obvious that the representation $\pi$ is bounded on $\CH$,
$\CU_q(su(2))$-equivariant
and $[\pi(\a),\gamma] = 0$ for all $\a \in \CA(S^2_q)$.

The next ingredient we search for is the reality structure
given by an antilinear operator $J$. Denote by $\ket{l,m}_\pm$
the orthonormal basis in the first (resp. second) copy
of $\CH_{\oh}$ in $\CH$. We claim the following:
\th{Lemma 6.}{}{
For all real $p>0$ the operator:
\begin{equation}
J \ket{l,m}_\pm = i^{2m} p^{m} \ket{l,-m}_\mp, \label{defJ}
\end{equation}
satisfies $J^2=-1$, $\gamma J = - J \gamma$
and
\begin{equation}
\label{tocomm}
J\pi(\a)J\pi(\b) = \pi(\b)J\pi(\a)J, \;\;\; \forall \a, \b \in \CA(S^2_q).
\end{equation}
}
\Proof Since $m$ is fractional (hence $2m$ is odd) and $J$ is antilinear, we have:
\begin{eqnarray*}
J^2 \ket{l,m}_\pm &=& J \left( i^{2m} p^m \ket{l,-m}_\mp \right) \\
&=& -\, i^{2m} p^m J \ket{l,-m}_\mp =
-\, i^{2m} p^m i^{-2m} p^{-m} \ket{l,m}_\pm
= - \ket{l,m}_\pm.
\end{eqnarray*}
The second property is obvious from the construction.
The only nontrivial part is the verification of (\ref{tocomm}).
We have verified this explicitly
for two pairs $(\a, \b) = (B, B)$ and $(B, B^*)$ of generators,
which clearly suffices in view of the relation (\ref{combi}).
Also, as an independent check, we have verified (\ref{tocomm})
with the help of computer symbolic computation program
for all the generators of the algebra.
\sq

\th{Remark 7.}{}{Note that $J$ is not unitary unless $p=1$,
in which case it has the same form as
the classical charge conjugation on the sphere.
However, only for $p=q$, $J$ is $\CU_q(su(2))$-equivariant,
that is for all $h \in \CU_q(su(2))$:
\begin{equation}
h J = J (Sh)^*.
\end{equation}
In fact, it is the equivariance property [15,12], which
led us to this form of $J$.
}

Finally, we shall look for all possible Dirac operators $D$,
which satisfy the following properties:
\begin{itemize}
\item $D$ is (in general unbounded) selfadjoint operator with compact resolvent
defined on a dense domain in $\CH$ which contains $V\oplus V$,
\item $D$ is $\CU_q(su(2))$-equivariant, that is
$$ Dh = hD, \;\;\;\; \forall h \in \CU_q(su(2)),$$
\item $D\gamma = -\gamma D$,
\item the commutators $[D, \pi(\a)]$ are bounded for all $\a \in \CA(S^2_q)$,
\item $D J = J D$,
\item $D$ satisfies the first order condition:
\end{itemize}
\begin{equation}
\label{FOC}
\left[ [D, \pi(\a)], J\pi(\b)J \right] =
0\;, ~~\forall \a, \b \in \CA(S^2_q).
\end{equation}
~\\
The following is the main result of the paper.
\th{Theorem 8.}{}{All operators satisfying the conditions listed above
are of the following form:
\begin{equation}
\label{dirac}
~~~~D = \left(\! \begin{array}{ll} 0 & {\bar z}D_- \\ z D_+ & 0 \end{array}\!\! \right),
\end{equation}
{\it ~where}
$$ D_{\pm} \ket{l,m}_\pm =  [l +\frac{1}{2}] ~\ket{l,m}_\mp $$
{\it and~} $z\in {\Bbb C}\setminus \{ 0\}$.
}
~\\
\Proof First, from the commutation with $\gamma$ we
obtain that $D: \CH_\pm \to \CH_\mp$.
Then the equivariance enforces that:
\begin{equation}
D \ket{l,m}_\pm =   d_l^\pm \ket{l,m}_\mp,
\label{Dei}
\end{equation}
where $d_l^\pm$ are constants.
In addition, the condition for $D$ to be hermitian yields
\begin{equation}
d_l^\pm = \overline{d_l^\mp} .
\label{Deig}
\end{equation}
Now, using the form of $J$ from (\ref{defJ}) we
calculate $DJ-JD$:
\begin{eqnarray*}
(DJ - JD) \ket{l,m}_\pm &=&
D i^{2m} p^m \ket{l,-m}_\mp - J d_l^\pm \ket{l,m}_\mp  \\
&=& i^{2m} p^m (d_l^\mp - \overline{d_l^\pm}) \ket{l,m}_\pm = 0.
\end{eqnarray*}
Here we used (\ref{Deig}) in the final step.

Next, as a consequence of the derivation property of commutators
and of (\ref{tocomm}),
the first order condition is satisfied
for all  $\a, \b \in \CA(S^2_q)$ iff it is satisfied
for all pairs of generators.
Furthermore, in view of (\ref{combi}),
the Jacobi identity and skewsymmetry of commutators,
it suffices to verify (\ref{FOC}) just for
two pairs of generators $(\a, \b) = (B, B)$ and $(B, B^*)$.
This yields first the following recurrence relation for $d_l^+$:
\begin{equation}
d_l^+ + d_{l+2}^+ = [2] d_{l+1}^+.
\end{equation}
Hence,
\begin{equation}
d_l^+ = x q^l + y q^{-l},
\end{equation}
where $x$ and $y$ are arbitrary constants.
Next, (\ref{FOC}) also yields $x = - q y$,
thus giving us the form of the Dirac operator as in (\ref{dirac}),
where $z$ is a complex parameter
(we can rescale all $d_l^+$ by the same $z$ and all
$d_l^-$ by $\bar z$).
Using this form of $D$
we verify explicitly (again with the help of symbolic computation program)
that the first-order condition is satisfied
for all pairs of generators of $\CA(S^2_q)$.

The parameter $z$ must be nonzero in order for the resolvent of $D$ to be compact.
It is also clear that the only accumulation points of the spectrum
$\{\pm |z|[l+1/2]\}$ of $D$ are $\pm \infty$
which guarantees the compactness of the resolvent.
Finally, the Dirac operator $D$ is selfadjoint on the dense domain,
$$
\left\{ \psi =\sum_{l, m, \sigma} c_{l, m, \sigma}  \ket{l,m}_\sigma ~:~
c_{l, m, \sigma}\in {\Bbb C}\ ,~
\sum_{l, m, \sigma} (1+ |z|^2 [l+1/2]^2) |c_{l, m, \sigma}|^2  < \infty
\right\}\ ,
$$
where the sum runs over $l=\frac{1}{2}, \frac{3}{2}, \ldots$,
$m =-l,\ldots,l$ and $\sigma=\pm$.

For this form of $D$ we show that
the commutators $[D, \pi(\a)]$,
for the elements $\a$ of the algebra $\CA(S^2_q)$ are bounded.
To see this
it suffices to check it only for the generators.
We start with the generator $B$.
For a given $l$, there are only three nonzero matrix elements
$ {}_-\langle l+j, m+1| [D,B] \ket{l,m}_+ $ corresponding to $j=-1,0,1$.
We verify that they are bounded.
In the case $j=1$ we estimate:
\begin{equation}
\begin{array}{l}
\phantom{xxxx} |{}_-\langle l+1, m+1| [D, \pi(B)] \ket{l,m}_+| \\
~\\
\phantom{xxxx} = q^m
\frac{\sqrt{[l+m+1][l+m+2]}}{\sqrt{[2l+2]([4l+4]+[2][2l+2])}}
q^{-l-2} \left( [l + \frac{3}{2}] - q [l + \frac{1}{2}] \right) \\
~\\
\phantom{xxxx} \leq C_+ q^{m} q^{-l-m} q^{3l} q^{-l} q^{-l} = C_+,
\end{array}
\end{equation}
where $C_+$ is a constant.
Similarly,
for $j=0$:
\begin{equation}
\begin{array}{l}
\phantom{xxxx} |{}_-\langle l, m+1| [D, \pi(B)] \ket{l,m}_+| \\
~\\
\phantom{xxxx} = q^m \frac{\sqrt{[l+m+1][l-m]}}{\sqrt{q} [2l] [2l+2]}
\left( q^{-1} - q \right) [l + \frac{1}{2}] \\
~\\
\phantom{xxxx} \leq C_0 q^{m} q^{-l} q^{4l} q^{-l} \leq C_0 ,
\end{array}
\end{equation}
for a suitable constant $C_0$.
Finally,
for $j=-1$:
\begin{equation}
\begin{array}{l}
\phantom{xxxx} |{}_-\langle l-1, m+1| [D, \pi(B)] \ket{l,m}_+| \\
~\\
\phantom{xxxx} = q^m
\frac{\sqrt{[l-m][l-m-1]}}{\sqrt{[2l+2]([4l+4]+[2][2l+2])}}
q^{-l-2} q^{2l}  \left( [l + \frac{3}{2}] - q [l + \frac{1}{2}] \right) \\
~\\
\phantom{xxxx} \leq C_- q^{m} q^{-l+m} q^{3l} q^{-l} q^{-l} q^{2l} \leq C_-,
\end{array}
\end{equation}
where in the last line we have used $q^{-m} \leq q^{-l}$.
It can be verified directly that similar bounds hold also
for the generators $A$ and $B^*$, but this in fact
is a consequence of (\ref{combi}) and the fact that $\pi$
is a bounded $*$-representation.
\sq

Note that we get an excellent agreement with the picture of the classical spectral
triple on the commutative sphere. Using the notation of [2] (pp.407--419)
we can (at $q=1$) identify $ \ket{l,m}_\pm$ with $Y^\pm_{l,m}$. Similarly,
for $p=1$ our $J$ has the same form as the charge conjugation $C$ on spinors and one
finds the eigenvectors and eigenvalues of $D$ to have the correct $q=1$ limit.

\section*{5. Conclusions.}

It is quite surprising that despite the earlier mentioned common belief that
Connes' approach to noncommutative geometry does not cover quantum-group examples,
we were able to produce an example satisfying the algebraic axioms
of real spectral triples (noncommutative manifolds).
We were informed by Alain Connes that also the dimension axiom can be satisfied
(related to the asymptotic behavior of the eigenvalues of the Dirac operator).
We believe that elaborating on this example will prove helpful to reconcile
the noncommutative geometry with the $q$-geometry.

It should be noted that, unlike in most of other approaches [1,10,11,13],
we do not assume at the beginning any form of the Dirac operator
(apart from the requirement of equivariance) and we derive it step by step.
Although its eigenvalues are the same as, e.g., those in [1],
the essential difference is the explicit realization of all the data
on the Hilbert space.
This includes the representation of the algebra and the proof
that the resulting differential forms have bounded representation $[D, \pi(\a)]$,
as well as a presentation of reality and chirality operations.
Also, in comparison to [7,3] (apart from using distinct starting algebras)
our Hilbert space representation and the equivariance requirements are fundamentally
different. Moreover, an important new ingredient is the reality operator.

There are still important problems to solve, which we shall address elsewhere:
the analytic properties of the spectral triple (summability, regularity and finiteness)
as well as the existence of the Hochschild cycle giving the volume form
and Poincar\'e duality.
(For details on these notions and their applications see [2], 10.5.)
Last not least, we plan to study the details of the classical limit $q\to 1$,
the relation of $\CH_\pm$ and $\CH$ to q-deformed Hopf bundles [8],
and the equivariant Fredholm module obtained from our construction.
Also the systematic analysis of other Podle\'s spheres
(in particular, the equatorial quantum sphere) is in preparation [6].

\vspace{5mm}
{\small \noindent{\bf Acknowledgments:} The authors would like
to thank P.~M.~Hajac for discussions and comments on the ma\-nu\-script,
and G.~Landi for discussions.
The first author acknowledges the financial support of the `Geometric Analysis'
Research Training Network HPRN-CT-1999-00118 of the E.~C.
The second author would like to thank ESF for a travel grant
and SISSA for hospitality.}
\vspace{5mm}
\vspace{5mm}
\begin{center}
\bf References
\end{center}

\begin{itemize}
\item[{[1]}]
{\sc P.~N.~Bibikov, P.~P.~Kulish},
{\it Dirac operators on the quantum group ${\rm SU}\sb q(2)$ and
the quantum sphere},  (Russian. English, Russian summary)
Zap. Nauchn. Sem. St. Petersburg. Otdel. Mat. Inst. Steklov.\  (POMI)  245 (1997),
Vopr. Kvant. Teor. Polya i Stat. Fiz.\  14, 49--65,
283;  translation in  J.~Math.~Sci.\ (New York)  100  (2000) 2039--2050.

\item[{[2]}]
{\sc J. M. Gracia-Bond\'{\i}a, J. C. V\'arilly, H. Figueroa},
{\it Elements of Noncommutative Geometry},
Birkh\"auser Advanced Texts, Birkh\"auser, Boston, MA, 2001.

\item[{[3]}]
{\sc P.~S.~Chakraborty, A.~Pal}
{\it Equivariant spectral triples on the quantum $SU(2)$ group},
~math.KT/0201004.

\item[{[4]}]
{\sc A.~Connes}, {\it Noncommutative Geometry}, Academic Press, 1994.

\item[{[5]}]
{\sc A.~Connes},
{\it Noncommutative geometry and reality},
J.~Math.~Phys. 36,  (1995)  6194--6231.

\item[{[6]}]
{\sc L.~D\c abrowski, G.~Landi, M.~Paschke, A.~Sitarz},
{\it in preparation}\/,

\item[{[7]}]
{\sc D.~Goswami},
{\it Some Noncommutative Geometric Aspects of $SU_q(2)$},
~math-ph/0108003.

\item[{[8]}]
{\sc P.~M.~Hajac, S.~Majid},
{\it Projective module description of the q-monopole}
Commun. Math. Phys. 206, (1999) 247-264.

\item[{[9]}]
{\sc A.~Klimyk, K.~Schmuedgen},
{\it Quantum Groups and their Representations}, Springer,
New York, 1998.

\item[{[10]}]
{\sc K.~Ohta, H.~Suzuki},
{\it Dirac Operators on Quantum Two Spheres},
Mod. Phys. Lett. A9, (1994) 2325--2334.

\item[{[11]}]
{\sc R.~Owczarek},
{\it Dirac operator on the Podles sphere},
Int.~J.~Theor.~Phys. 40, (2001) 163--170.

\item[{[12]}]
{\sc M.~Paschke},
{\em \"Uber nichtkommutative Geometrien, ihre Symmetrien und
etwas Hochenergiephysik}, Ph.D.~Thesis, Mainz 2001.

\item[{[13]}]
{\sc A.~Pinzul, A.~Stern},
{\it Dirac Operator on the Quantum Sphere},
Phys.~Lett.~B  512 (2001) 217--224.

\item[{[14]}]
{\sc P.~Podle\'{s}},
{\it Quantum Spheres},
Lett.~Math.~Phys. 14 (1987) 521--531.

\item[{[15]}]
{\sc A.~Sitarz},
{\it Equivariant Spectral Triples},
{in this volume}.

\item[{[16]}]
{\sc S.~L.~Woronowicz},
{\it Twisted $SU(2)$ group. An example of a
noncommutative differential calculus},
Publ.\ RIMS, Kyoto University, 23 (1987), 117--181.
\end{itemize}
\end{document}